\documentclass[12pt]{amsart}
\usepackage{amssymb,amsthm,rotating, mathtools}
\usepackage[all]{xy}
\numberwithin{equation}{section}
 \newtheorem{thm}[equation]{Theorem}
\newtheorem{theorem}[equation]{Theorem}
 \newtheorem{defn}[equation]{Definition}
 \newtheorem{prop}[equation]{Proposition}

\newtheorem{cor}[equation]{Corollary}

 \theoremstyle{definition}

 \newtheorem{remark}[equation]{Remark}

 \newtheorem{example}[equation]{Example}

\DeclareMathOperator{\Ext}{Ext}
\DeclareMathOperator{\Hom}{Hom}

\DeclareMathOperator{\GL}{GL}
\DeclareMathOperator{\Gl}{GL}

\DeclareMathOperator{\sgn}{sgn}

\DeclareMathOperator{\Sym}{Sym}
\DeclareMathOperator{\Tor}{Tor}

\DeclareMathOperator{\diag}{diag}

\newcommand{\DOT}{\setlength{\unitlength}{1pt}\begin{picture}(2.5,2)
                  (1,1)\put(2,3.5){\circle*{2}}\end{picture}}

\newcommand{\del}{\partial}
\newcommand{\HH}{{\rm HH}}
\newcommand{\HHD}{\HH^{\DOT}}
\newcommand{\Wedge}{\textstyle\bigwedge}
\newcommand{\CC}{\mathbb{C}}
\newcommand{\NN}{\mathbb{N}}
\newcommand{\ZZ}{\mathbb{Z}}
\newcommand{\ot}{\otimes}

\newcommand{\cohg}{{\rm HH}^{\DOT}\bigl(S(V),S(V)\bar{g}\bigr)}
\newcommand{\ep}{\epsilon}
\newcommand{\scoochmore}{\hspace{-.5ex}}

\newcommand{\ignore}[1]{\relax}

\newcommand{\ta}{\Upsilon}
\title[Chain maps]
{Quantum differentiation and chain maps of bimodule complexes
}
\date{July 15, 2009}
\author{Anne V.\ Shepler}
\address{Department of Mathematics, University of North Texas,
Denton, Texas 76203, USA}
\email{ashepler@unt.edu}
\author{Sarah Witherspoon}
\address{Department of Mathematics\\Texas A\&M University\\
College Station, Texas 77843, USA}\email{sjw@math.tamu.edu}
\thanks{The first author was partially supported by NSF grants
\#DMS-0402819 and \#DMS-0800951 and a research fellowship from the 
Alexander von Humboldt Foundation. 
The second author was partially supported by NSA grant
\#H98230-07-1-0038 and NSF grant
\#DMS-0800832.
Both authors were jointly supported by Advanced Research Program Grant 010366-0046-2007
from the Texas Higher Education Coordinating Board.}

\begin{document}
\maketitle

\begin{abstract}
We consider a finite group acting on a vector space and
the corresponding skew group algebra 
generated by the group and the symmetric algebra of the
space.  This skew group algebra illuminates the resulting orbifold and serves as a replacement 
for the ring of invariant polynomials, especially in the eyes of cohomology.  
One analyzes the Hochschild cohomology of the skew group algebra 
using isomorphisms which convert between resolutions.
We present an explicit chain map from the bar resolution to the Koszul resolution of the
symmetric algebra which induces various isomorphisms on Hochschild homology
and cohomology, some of which have appeared in the literature before.
This approach unifies previous results on homology and
cohomology of both the symmetric algebra and skew group algebra.
We determine induced
combinatorial cochain maps which invoke quantum differentiation
(expressed by Demazure-BBG operators).
\end{abstract}

%%%%%%%%%%%%%%%%%%%%%%%%%%%%%%%%%%%%%%%%%%%%%%%%%%%%%%%%%%%%%%
%%%%%%%%%%%%%%%%%%%%%%%%%%%% SECTION %%%%%%%%%%%%%%%%%%%%%%%%%
%%%%%%%%%%%%%%%%%%%%%%%%%%%%%%%%%%%%%%%%%%%%%%%%%%%%%%%%%%%%%%
\section{Introduction}

Let $G$ be a finite group acting linearly on a finite dimensional, complex vector space
$V$.
The skew group algebra $S(V)\# G$ is a natural semi-direct 
product of $G$ with the symmetric algebra $S(V)$ (a polynomial ring).
It serves as a valuable, albeit noncommutative, 
replacement for the invariant ring $S(V)^G$ in geometric settings, 
as it encodes the abstract group structure of $G$ as well as its action on $V$.
The cohomology of $S(V)\# G$
informs various areas of mathematics (for example,
geometry, combinatorics, representation theory, 
and mathematical physics).  In particular, the Hochschild cohomology
of $S(V)\#G$ governs its deformations, which include graded Hecke algebras,
symplectic reflection algebras, and Cherednik algebras.

The orbifold $V/G$ may be realized as an algebraic variety
whose coordinate ring is the ring of invariant functions $S(V^*)^G$ 
on the dual space $V^*$, which is the center of $S(V^*)\# G$
when $G$ acts faithfully. 
(For details, see Harris~\cite{Harris}.)  
The variety $V/G$ is nonsingular exactly when the action of $G$ on $V$
is generated by reflections.
Geometers and physicists are interested in resolving the singularities
of $V/G$ with a smooth variety $X$ and examining the coordinate ring of $X$ instead
of $S(V)^G$.  In situations they study, the skew group algebra
$S(V)\#G$ serves as a replacement for the coordinate ring of $X$; indeed, 
Hochschild cohomology sees no difference between these rings
(see C\u{a}ld\u{a}raru, Giaquinto, 
and the second author~\cite{CaldararuGiaquintoWitherspoon}).
Connections with representation theory are still unfolding, for example, see
Gordon and Smith~\cite{GordonSmith}.

The Hochschild cohomology $\HH^{\DOT}(A)$ of any algebra $A$ over a field $k$
is the space $\Ext^{\DOT}_{A\ot A^{op}}(A,A)$. 
The cup product and Gerstenhaber bracket on Hochschild cohomology
are both defined initially on the bar resolution, a natural $A\ot A^{op}$-free
resolution of $A$. The cup product has another description as Yoneda
composition of extensions of modules, which can be transported to any
other projective resolution.  However, the Gerstenhaber bracket has 
resisted such a general description. Instead, one commonly computes
$\HH^{\DOT}(A)$ using a more convenient resolution, then one finds and uses
relevant chain maps to lift the Gerstenhaber bracket from the bar resolution.
The case $A=S(V)\# G$ is complicated further because one 
does not work with resolutions of $A$ directly, but instead one derives information 
from resolutions of the symmetric algebra $S(V)$.

In this paper, we begin this task by constructing explicit chain maps 
which encode traffic between resolutions used to describe the Hochschild
cohomology of $A=S(V)\# G$.  
Our maps convert between the bar and
Koszul resolutions of the polynomial ring $S(V)$, and serve as a tool for investigating
the homology and cohomology of $S(V)$ with coefficients in any bimodule.
Specifically, 
the Koszul resolution of the polynomial ring $S(V)$ embeds naturally into the bar resolution.
We define an explicit chain map, depending on a choice of basis, 
giving a quasi-inverse to this embedding.
We study in particular the induced maps on the Hochschild cohomology
$\HH^{\DOT}(S(V),S(V)\# G)$.
We give an  elegant, combinatorial description of the induced map
on cochains in terms of scaled Demazure
(BBG) operators (or quantum partial differential operators, see
Definition~\ref{tau}). 
We describe the induced maps
on Hochschild homology as well.
(These combinatorial descriptions are useful for computations,
which we pursue in other articles.)  
The cohomology 
$\HHD(S(V)\# G)$ manifests as the $G$-invariant subspace of
$\HHD(S(V),S(V)\# G)$
in characteristic 0.
We thus obtain isomorphisms of homology
and cohomology that allow one to transfer structures defined on the bar
resolution to the complexes standardly used to describe
$\HHD(S(V)\# G)$.

In Section~\ref{Prelim}, we establish notation and deploy
the Hochschild cohomology $\HH^{\DOT}(S(V)\# G)$ 
in terms of both the Koszul and bar resolutions of $S(V)$.  
We introduce a
combinatorial map $\ta$ on cochains in Section~\ref{Quantum}.
This combinatorial converter $\ta$  takes vector forms (tagged by
group elements) to twisted quantum differential operators.
In Section~\ref{chainmap},
we give a technical formula for explicit chain maps from the
bar resolution to the Koszul resolution (Definition~\ref{psi}), 
which is valid over an arbitrary ground field.
These specific chain maps each induce an inverse to the embedding 
of the Koszul resolution into the bar resolution after taking homology or cohomology.
(Indeed, after applying functors $\ot$ or $\Hom$, we recover some
chain maps given in the literature 
for converting between complexes expressing Hochschild homology
and cohomology---see Section~\ref{intheliterature}.)
In Section~\ref{merits},
we deduce that our combinatorial converter $\ta$ defines automorphisms of
cohomology
by showing that it is induced by the chain maps of Section~\ref{chainmap}.
We present similar automorphisms of homology 
(using quantum differentiation) in Section~\ref{intheliterature}.

Our approach presents an immediate and obvious advantage:
We define one primitive map between resolutions and then apply various
functors that automatically give (co)chain maps in a variety of settings.
We do not need to give separate proofs (depending on context)
showing that these induced maps are chain maps,
as such results follow immediately from the general theory.
This uniform treatment provides a clear channel
for navigating 
between chain and cochain complexes.  
Indeed, we use this channel in~\cite{SWring}
to explore the algebraic structure of 
$\HH^{\DOT}(S(V)\#G)$ under the cup product.
We shall use this approach in future work to obtain results on the Gerstenhaber bracket.

Some results in this paper are valid over a field of arbitrary 
characteristic, while others assume the ground field is the complex numbers, $\CC$.
We have tried to state carefully requirements on the field throughout.
The reader should note that whenever we work over $\CC$, we could instead work
over any field containing the eigenvalues of the action of $G$ on $V$ in which
$|G|$ is invertible.
All tensor and exterior products will be taken over the ground
field unless otherwise indicated.

%%%%%%%%%%%%%%%%%%%%%%%%%%%%%%%%%%%%%%%%%%%%%%%%%%%%%%%%%%%%%%
%%%%%%%%%%%%%%%%%%%%%%%%%%%% SECTION %%%%%%%%%%%%%%%%%%%%%%%%%
%%%%%%%%%%%%%%%%%%%%%%%%%%%%%%%%%%%%%%%%%%%%%%%%%%%%%%%%%%%%%%
\section{Preliminary material}
\label{Prelim}

In this section, 
we work over the complex numbers $\CC$,
although the definitions below of Hochschild cohomology, bar resolution,
and Koszul resolution are valid over any ground field.

Let $G$ be a finite group and $V$ a (not necessarily faithful) $\CC G$-module.
Let ${}^g v$ denote the image of $v\in V$ under the action of 
$g\in G$.  We work with the induced group action on all maps
throughout this article:
For any map $\theta$ and element $h\in \Gl(V)$, 
we define the map $ {}^h\theta $ 
by $(^h\theta)(v) := \hphantom{\,}^h(\theta(^{h^{-1}}v))$ for all $v$.
Let $V^*$ denote the vector space dual to $V$ with
the contragredient (i.e., dual) representation.
For any basis $v_1,\ldots,v_n$ of $V$, let
$v_1^*,\ldots,v_n^*$ be the dual basis of $V^*$.
Let $V^G=\{v \in V: {}^{g}v=v \text{ for all } g\in G\}$, the set of 
$G$-invariants in $V$. 
For any $g\in G$, let $Z(g)=\{h\in G: gh=hg\}$, 
the centralizer of $g$ in $G$,
and let $V^g = \{v\in V: {}^gv=v\}$, the $g$-invariant subspace of $V$.

The {\bf skew group algebra} $S(V)\# G$ is the vector space
$S(V)\ot \CC G$ with multiplication given by
$$
  (a\ot g)(b\ot h)=a ( {}^{g}b) \ot gh
$$
for all $a,b\in S(V)$ and $g,h\in G$.
We abbreviate $a\ot g$ by $a\overline{g}$ ($a\in S(V)$, $g\in G$)
and $a\ot 1, \ 1\ot g$ simply by $a, \ \overline{g}$, respectively.
An element $g\in G$ acts on $S(V)$ by an inner automorphism in $S(V)\# G$: 
$ \ \overline{g} a (\overline{g})^{-1} = ( {}^ga)\overline{g} (\overline{g})
^{-1} = {}^ga$ for all $a\in A$.

The {\bf Hochschild cohomology} of a $\CC$-algebra $A$ (such as $A=S(V)\# G$),
with coefficients in an $A$-bimodule $M$,
is the graded vector space
$\HH^{\DOT}(A,M)=\Ext^{\DOT}_{A^e}(A,M)$, where $A^e=A\ot A^{op}$ acts
on $A$ by left and right multiplication.
This cohomology may be expressed in terms of the 
{\bf bar resolution}, the following free $A^e$-resolution of $A$:
\begin{equation}\label{barcomplex}
  \cdots\stackrel{\delta_3}{\longrightarrow} A^{\ot 4}
\stackrel{\delta_2}{\longrightarrow} A^{\ot 3}
\stackrel{\delta_1}{\longrightarrow} A^e
\stackrel{m}{\longrightarrow} A \rightarrow 0,
\end{equation}
where $\delta_p(a_0\ot\cdots\ot a_{p+1}) = \sum_{j=0}^p (-1)^j a_0
\ot\cdots\ot a_ja_{j+1}\ot\cdots\ot a_{p+1}$, and $\delta_0=m$ is
multiplication.
We apply $\Hom_{A^e}( - , M)$ to obtain a cochain complex whose homology is
$\HH^{\DOT}(A,M)$.  
If $M=A$, we abbreviate $\HH^{\DOT}(A)=\HH^{\DOT}(A,A)$.
For each $p$, $\Hom_{A^e}(A^{\ot (p+2)},A)\cong \Hom_{\CC}(A^{\ot p}, A)$,
and we identify these two spaces of $p$-cochains
throughout this article.
The graded vector space $\HH^{\DOT}(A)$ admits both a cup product and a
graded Lie bracket under which it becomes a Gerstenhaber algebra.
In this article, we develop automorphisms of cohomology converting between 
resolutions.  These automorphisms will be used in subsequent publications
to explore the algebraic structure 
of $\HH^{\DOT}(S(V)\# G)$ under these two operations.

%%%%%%%%%%%%%%%%%%%%%%%%%%%%%%%%%%%%%%%%%%%%%%%%%555
%%%%%%%%%%%%%%%%%%%%%%%%%%%%%%%%%%%%%%%%%%%%%%%%%%%%%5
%%%%%%%%%%%%%%%%%%%%%%%%%%%%%%%%%%%%%%%%%%%%%%%%%%%%%%%%%5
\subsection*{Hochschild cohomology of $S(V)\#G$}
\label{Isomorphisms}
The graded vector space structure of
$\HH^{\DOT}(S(V)\# G)$ was determined by 
Farinati~\cite{Farinati} and Ginzburg and Kaledin~\cite{GinzburgKaledin} 
when $G$ acts faithfully on $V$.
The same techniques apply to nonfaithful actions.
The following statements are valid only when the characteristic does not divide
the order of $G$.
(Otherwise, the cohomology is more complicated
as the group algebra of $G$ may itself not be semisimple.)
Let $\mathcal C$ be  a set of representatives of the conjugacy classes
of $G$. 
A consequence of a result of 
\c{S}tefan~\cite[Cor.\ 3.4]{Stefan} posits a natural $G$-action
giving
the first of a series of 
isomorphisms of graded vector spaces:
\begin{equation}\label{isos}
\begin{aligned}
  \HH^{\DOT}(S(V)\# G) \ \ 
&\cong \ \ 
\HH^{\DOT}(S(V), S(V)\#G)^G\\
\ \ &\cong \ \ 
 \Biggl(\bigoplus_{g\in G} \HH^{\DOT}(S(V), S(V)\bar{g})
\Biggr)^{G}\\
\ \ &\cong \ \ 
\bigoplus_{g\in{\mathcal C}}\cohg^{Z(g)}.
\end{aligned}
\end{equation}
Specifically,
the action of $G$ on $V$ extends naturally to the bar complex of $S(V)$
and commutes with the differentials,
and thus it induces a natural action on Hochschild cohomology $\HH^{\DOT}(S(V),S(V)\#G)$
(for which we also use the action of $G$ on $S(V)\# G$ by inner automorphisms).
The subspace of $G$-invariants of this action is denoted by 
$\HH^{\DOT}(S(V),S(V)\#G)^G$.  
(Equivalently, this may also be defined via any other
choice of $G$-compatible resolution used to compute cohomology, 
see~\cite[Section 2]{Stefan}, for example).

The second isomorphism of~(\ref{isos}) surfaces simply because 
the $S(V)^e$-module $S(V)\# G$ decomposes into the direct sum of 
$S(V)^e$-modules $S(V)\overline{g}$, and cohomology
preserves direct sums.  (The isomorphism arises of course at the
cochain level, as the $\Hom$-functor preserves direct sums.)  We identify
$\HH^{\DOT}(S(V), S(V)\#G)$ with
$\bigoplus_{g\in G} \HH^{\DOT}(S(V), S(V)\bar{g})$ when convenient throughout
this article.
Note that $G$ permutes the components in the direct sum in accordance
with the conjugation action of $G$ on itself.
Thus for each $g\in G$, the subgroup $Z(g)$ fixes the $g$-component 
$\HH^{\DOT}(S(V),S(V)\overline{g})$ setwise. 
The third isomorphism of~(\ref{isos}) canonically projects 
onto a set of representative summands.

One may use the Koszul resolution for $S(V)$ to determine each $g$-component
$\HH^{\DOT}(S(V), S(V)\overline{g})$ in the last line of~(\ref{isos}) above.
The {\bf Koszul resolution}, denoted $K_{\DOT}(S(V))$, is given by 
$K_0(S(V))=S(V)^e$, $K_1(S(V))= S(V)^e\ot V$, and for each $p\geq 2$,
\begin{equation}\label{koszul-res1}
  K_p (S(V)) = \bigcap_{j=0}^{p-2} S(V)^e\ot (V^{\ot j}\ot R\ot V^{\ot (p-j-2)}),
\end{equation}
where $R$ is the subspace of $V\ot V$ spanned by all $v\ot w - w\ot v$
($v,w\in V$) (e.g., see Braverman and Gaitsgory~\cite{BravermanGaitsgory}).
This is a subcomplex of the bar resolution~(\ref{barcomplex}) for $S(V)$.
For any choice of basis $v_1,\ldots,v_n$ of $V$, it is equivalent to
the Koszul resolution corresponding to the regular sequence
$\{v_i\ot 1-1\ot v_i\}^n_{i=1}$ in $S(V)^e$:
\begin{equation}\label{koszul-res2}
  K_{p}(\{v_i\ot 1 - 1\ot v_i\}^n_{i=1}) \cong S(V)^e \ot \Wedge^{p}(V),
\end{equation}
a free $S(V)^e$ resolution of $S(V)$ (e.g., see Weibel~\cite[\S4.5]{Weibel}).
The differentials are given by 
\begin{equation}\label{koszul-diff}
d_p(1\ot 1\ot v_{j_1}\wedge\cdots\wedge v_{j_p}) = 
   \sum_{i=1}^p (-1)^{i+1} (v_{j_i}\ot 1 - 1\ot v_{j_i})\ot
   (v_{j_1}\wedge\cdots\wedge \hat{v}_{j_i}\wedge\cdots\wedge v_{j_p}).
\end{equation}
The canonical inclusion of the Koszul resolution~(\ref{koszul-res1})
into the bar resolution~(\ref{barcomplex}) for $S(V)$ is then given on 
resolution~(\ref{koszul-res2})
by the chain map $$\Phi : S(V)^e\ot \Wedge^{\DOT}(V)\rightarrow
S(V)^{\ot (\DOT +2)}\ ,$$ defined by 
\begin{equation}\label{phik}
  \Phi_p(1\ot 1\ot v_{j_1}\wedge\cdots\wedge v_{j_p})
   = \sum_{\pi\in \Sym_p}\sgn(\pi)\ot v_{j_{\pi(1)}}\ot\cdots\ot
    v_{j_{\pi(p)}}\ot 1
\end{equation}
for all $v_{j_1},\ldots,v_{j_p}\in V$, $p\geq 1$,
where $\Sym_p$ denotes the symmetric group
on $p$ symbols.
Note that by its definition, $\Phi$ is invariant
under the action of $\GL(V)$, i.e.,\ $^h\Phi = \Phi$ for all $h$ in $\Gl(V)$.

Any chain map
   $\Psi_p  :  S(V)^{\ot (p+2)}\rightarrow S(V)^e\ot \Wedge^p V$
from the bar resolution to the Koszul resolution yields
a commutative diagram:
$$
%\begin{xy}
\xymatrix{
\cdots \ar[r] & S(V)^{\ot 4}\ar[r]^{\delta_2}\ar@<-2pt>[d]_{\Psi_2} 
               & S(V)^{\ot 3}\ar[r]^{\delta_1}\ar@<-2pt>[d]_{\Psi_1} 
               & S(V)^e \ar[r]^m \ar[d]_{=}
               & S(V) \ar[r] \ar[d]_{=} & 0\\
\cdots \ar[r] & S(V)^e\ot \Wedge ^2 V \ar[r]^{d_2}\ar@<-2pt>[u]_{\Phi_2}
               & S(V)^e\ot \Wedge ^1 V \ar[r]^{\hspace{.6cm}d_1}\ar@<-2pt>[u]_{\Phi_1}
               & S(V)^e \ar[r]^{m} \ar[u]
               & S(V)\ar[r] \ar[u] & 0 .
}
%\end{xy}
$$
(In Definition~\ref{psi} below, we explicitly 
define a map $\Psi$ depending on a choice of basis of $V$.)
Such maps $\Phi$ and $\Psi$ necessarily induce inverse isomorphisms on 
cohomology $\HH^{\DOT}(S(V),M)$ for any $S(V)$-bimodule $M$, upon
applying $\Hom_{S(V)^e}( - , M)$.
(Similarly for homology; see Section~\ref{intheliterature} below.)
After identifying
$ \Hom_{S(V)^e}(S(V)^{\ot (p+2)},M)$ with 
$\Hom_{\CC}(S(V)^{\ot p}, M)$ 
and
$\Hom_{S(V)^e}(S(V)^e \ot \Wedge^{p} V ,M)$ 
with $\Hom_{\CC}(\Wedge^{p} V, M)$ for all $p$, we obtain the following
commutative diagram:
\begin{equation}\label{big-diagram}
\xymatrix{
 \Hom_{\CC}(S(V)^{\ot p}, M) \ar[r]^{\delta_{p-1}^*}\ar@<2pt>[d]^{\Phi^*_p}
     & \Hom_{\CC}(S(V)^{\ot (p+1)}, M)\ar@<2pt>[d]^{\Phi^*_{p+1}}\\
  \Hom_{\CC}( \Wedge ^p V , M)\ar[r]^{d_{p-1}^*}  \ar@<2pt>[u]^{\Psi^*_p}
    & \Hom_{\CC}( \Wedge^{p+1} V, M)\ar@<2pt>[u]^{\Psi^*_{p+1}}\ \ \ .
}
\end{equation}
The maps $\Phi \Psi$ and $\Psi \Phi$ are each homotopic to an identity map
by the Comparison Theorem,
and thus $\Phi_p^*$ and $\Psi_p^*$ induce inverse automorphisms on  
the cohomology $\HH^{p}(S(V),M)$
(see the proof of Lemma 2.4.1 of Weibel~\cite{Weibel}).
In this paper, we primarily consider the $S(V)^e$-modules
$M=S(V)\# G$ and $M=S(V)\overline{g}$ for $g$ in $G$.

We transfer the map $\Phi$ and any chain map $\Psi$ to
the Hochschild cohomology of the full skew group algebra, $\HH^{\DOT}(S(V)\# G)$,
using the isomorphisms of~(\ref{isos}).
Set $M=S(V)\# G$ and let $\Phi^*$ and $\Psi^*$ denote the
induced maps on the cohomology
$$\HH^{\DOT}(S(V),S(V)\# G)\cong \bigoplus_{g\in G}\, \cohg \ .$$
For each $g$ in $G$, denote the restrictions to $\cohg $
by $\Phi_g^*$ and $\Psi_g^*$, respectively, 
so that
$$
\Phi^* = \bigoplus_{g\in G} \Phi_g^*
\quad\quad\text{and}\quad\quad
\Psi^* = \bigoplus_{g\in G} \Psi_g^*\ .
$$
The maps $\Phi^*$ and $\Psi^*$ behave nicely with respect to the action of $G$:
\vspace{2ex}
%%%%%%%%%%%%%%%%%%%%%%%%%%%%%%55
\begin{prop}
\label{psiInvariant}
Let $\Psi$ be any choice of chain map from the bar resolution~(\ref{barcomplex})
to the Koszul resolution~(\ref{koszul-res2}). 
The cochain maps $\Phi^*$ and $\Psi^*$
are inverse automorphisms on the cohomology $\HH^{\DOT}(S(V),S(V)\# G)$ 
converting between expressions arising 
from the Koszul resolution and from the bar resolution.
In addition,
\begin{itemize}
\item[(1)]
For any $g\in G$,
the maps $\Phi_g^*$ and $\Psi_g^*$
on the cohomology $\cohg$
are invariant under the centralizer $Z(g)$ of $g$ in $G$, and
the maps $\Phi^*$ and $\Psi^*$ on
$\bigoplus_{g\in G}\, \cohg$ are invariant under $G$;
\item[(2)]
The maps  $\Phi^*$ and $\Psi^*$ induce inverse automorphisms on 
 the graded vector space
$$
\hspace{.6cm}   \bigoplus_{g\in{\mathcal C}}\, \bigl(\cohg\bigr)^{Z(g)} \cong
    \left(\bigoplus_{g\in G}\HH^{\DOT}(S(V),S(V)\overline{g})\right)^G
   \cong \HH^{\DOT}(S(V)\# G).
$$
\end{itemize}
\end{prop}
%%%%%%%%%%%%%%%%%%%%%%%%%%%%%%%%%%%
\vspace{1ex}
\begin{proof}
As explained after Diagram~(\ref{big-diagram}),
the maps $\Phi_g^*$ and $\Psi_g^*$
are inverse isomorphisms on the cohomology $\cohg$.
By its definition, $\Phi$ is invariant under the action of $\GL(V)$, 
and so the map $\Phi^*$ on $\HH^{\DOT}(S(V),S(V)\#G)$
is invariant under $G$, and  
the map $\Phi_g^*$ on $\cohg$ is invariant under $Z(g)$. 
Fix some $h$ in $Z(g)$ and consider the map $^h(\Psi_g^*)$.
As maps on the cohomology $\cohg$, 
$$
1\ =\ ^h ( \Phi_g^*\ \Psi_g^* ) 
 \ =\ ^h(\Phi_g^*)\ ^h(\Psi_g^*)
 \ =\ \Phi_g^*\ ^h(\Psi_g^*) , 
$$
thus $^h(\Psi_g^*)$ is also inverse to $\Phi_g^*$.
Hence $^h(\Psi_g^*) = \Psi_g^*$ (since the inverse is unique) as maps on cohomology, 
for all $h$ in $Z(g)$, and $\Psi_g^*$ is also $Z(g)$-invariant. 
Thus statement (1) holds. 
As a consequence, we may restrict both $\Phi_g^*$ and $\Psi_g^*$ 
to the graded vector space  $\bigl(\cohg\bigr)^{Z(g)}$.
Applying the isomorphisms~(\ref{isos}), we obtain the statement (2). 
\end{proof}
\vspace{2ex}
%%%%%%%%%%%%%%%%%%%%%%%%%%%%%%%%%%%%%%%%%%%%%%%%%%%%%%%%%%%%

The cohomology $\HHD(S(V),S(V)\# G)$
arising from the Koszul resolution~(\ref{koszul-res2}) 
of $S(V)$ may be viewed as a set of
{\bf vector forms} on $V$ {\bf tagged} by group elements of $G$.
Indeed, we identify
$\Hom_{\CC}(\Wedge^p V, S(V)\overline{g})$ with 
$S(V)\overline{g}\ot \Wedge^p V^*$ for each $g$ in $G$,
and recognize 
the set of {\bf cochains} 
derived from the Koszul resolution as (see Diagram~\ref{taudiagram} below)
\begin{equation}\label{C-cochains}
C^{\DOT}:=\bigoplus_{g\in G} C_g^{\DOT},
\quad\text{where}\quad 
C_g^p :=  S(V)\bar{g} \otimes \Wedge^{p} V^*
\ .
\end{equation}

\newpage

%%%%%%%%%%%%%%%%%%%%%%%%%%%%%%%%%%%%%%%%%%%%%%%%%%%%%%%%%%%%%%
%%%%%%%%%%%%%%%%%%%%%%%%%%%% SECTION %%%%%%%%%%%%%%%%%%%%%%%%%
%%%%%%%%%%%%%%%%%%%%%%%%%%%%%%%%%%%%%%%%%%%%%%%%%%%%%%%%%%%%%%
\section{Quantum 
differentiation and a combinatorial converter  map}
\label{Quantum}

One generally uses the
Koszul resolution of $S(V)$ to compute Hochschild cohomology, 
but some of the algebraic structure 
of its cohomology is defined using the bar resolution instead.
We thus define automorphisms of cohomology which convert between
resolutions.
In Equation~(\ref{phik}), we defined the familar inclusion map $\Phi$ 
from the Koszul resolution to the bar resolution.
But in order to transfer algebraic structure, 
we need chain maps in both directions.
In Section~\ref{chainmap}, we shall construct
explicit chain maps $\Psi$ from the bar resolution to the Koszul resolution,
which will then induce cochain maps $\Psi^*$.
These maps $\Psi$ are somewhat unwieldy, however.
Thus, in this section, we first define a more 
elegant and handy  map $\ta$ on cochains 
using quantum differential operators (alternatively, Demazure operators).  
In Theorem~\ref{psi=tau}, we prove that $\ta = \Psi^*$ as maps on
cocycles, for our specific construction of a chain map $\Psi$ from the bar resolution 
to the Koszul resolution of $S(V)$.  
This impies that the map $\ta$ is itself a cochain map, and that
$\ta$ is in fact equal 
to $\Psi^*$ on cohomology,  
for {\em any} choice of chain map $\Psi$ from the bar resolution 
to the Koszul resolution of $S(V)$.
This development allows us to deduce important
properties of the expedient map $\ta$ (useful for computations)
from the elephantine map $\Psi^*$.
In this section, we work over the complex numbers $\CC$.

Given any basis $v_1, \ldots, v_n$ of $V$, and any 
complex number $\epsilon\neq 1$, we 
define the $\epsilon$-quantum partial
differential operator with respect to
$v:=v_i$ as the scaled Demazure
(BBG) operator 
$\del_{v, \epsilon}:S(V)\rightarrow S(V)$
given by
\begin{equation}\label{qpd}
\del_{v, \epsilon}(f)\ = \ (1-\epsilon)^{-1}\ \frac{f-\, ^s\hspace{-.5ex}f}{v}
\ = \ \frac{f-\, ^s\hspace{-.5ex}f}{v-\, ^sv}\ ,
\end{equation}
where $s\in\Gl(V)$ is the reflection whose matrix with respect to the 
basis $v_1, \ldots, v_n$ is 
$\diag(1,\ldots, 1, \epsilon, 1, \ldots, 1)$
with $\epsilon$ in the $i$th slot.
Set $\del_{v,\ep}=\del/\del v$, the usual partial differential operator
with respect to $v$, when $\ep=1$.

%%%%%%%%%%%%%%%%%%%%%%%%%%%%%%%%5
\vspace{3ex}
\begin{remark}
The quantum partial differential operator $\del_{v,\epsilon}$ 
above coincides with the usual definition of
quantum partial differentiation:
One takes the ordinary partial derivative with respect to $v$ but instead
of multiplying each monomial by its degree $k$ in $v$, one multiplies by
the quantum integer $[k]_{\epsilon}$, where
$$
   [k]_{\epsilon} := 1 + \epsilon + \epsilon^2 +\cdots + \epsilon^{k-1}.
$$
Let us check explicitly that these two definitions coincide. 
For $v=v_1$, $\ep\neq 1$,
$$
\begin{aligned}
\del_{v,\epsilon}(v_1^{k_1}v_2^{k_2}\cdots v_n^{k_n})
&=\frac{(v_1^{k_1}v_2^{k_2}\cdots v_n^{k_n})
  -\, ^s(v_1^{k_1}v_2^{k_2}\cdots v_n^{k_n})}{v_1-\, ^sv_1}\\
&=\frac{v_1^{k_1}v_2^{k_2}\cdots v_n^{k_n}
  -\epsilon^{k_1}v_1^{k_1}v_2^{k_2}\cdots v_n^{k_n}}{v_1-\epsilon v_1}\\
&=\frac{(1-\epsilon^{k_1}) v_1^{k_1}v_2^{k_2}\cdots v_n^{k_n}}
   {(1-\epsilon)v_1}\\ 
&=[k_1]_{\epsilon}\ \  v_1^{k_1-1}v_2^{k_2}\cdots v_n^{k_n}.
\end{aligned}
$$
\end{remark}
\vspace{3ex}
%%%%%%%%%%%%%%%%%%%%%%%%%%%%%%%%%%%%%%%%%%%%%%%%%%%%%%%%

We are now ready to construct the map $\ta$  
taking vector forms (tagged by group elements)
to twisted quantum differential operators.  We define
$\ta$ on cochains $C^{\DOT}$ (see~(\ref{C-cochains}))
so that the following diagram
commutes for $M=S(V)\# G$:
\begin{equation}\label{taudiagram}
\xymatrix{
 \Hom_{\CC}(S(V)^{\ot p}, M) \ar[rr]^{\delta^*_{p-1}\hspace{.2cm}}\ar@<2pt>[d]^{\Phi^*_p}
     && \Hom_{\CC}(S(V)^{\ot (p+1)}, M)\ar@<2pt>[d]^{\Phi^*_{p+1}}\\
  \displaystyle{\bigoplus_{g\in G}}\ S(V)\overline{g}\ot \Wedge^p V^*
\ar[rr]^{d^*_{p-1}\hspace{.2cm}}  \ar@<2pt>[u]^{\ta_p}
  && \ \ \displaystyle{\bigoplus_{g\in G}}\ S(V)\overline{g}\ot\Wedge^{p+1}V^*\ar@<2pt>[u]^{\ta_{p+1}}
\ \ \ .
}
\end{equation}
First, some notation.
For $g$ in $G$, fix a basis $B_g=\{v_1,\ldots, v_n\}$
of $V$ consisting of eigenvectors of $g$
with corresponding eigenvalues $\epsilon_1, \ldots, \epsilon_n$.
Decompose $g$ into a product of reflections diagonal in this basis: 
Let $g=s_1\cdots s_n$ where each $s_i$ is either the identity or a reflection 
defined by
${}^{s_i}v_j = v_j$ for $j\neq i$ and ${}^{s_i}v_i=\epsilon_i v_i$.
Let $\del_{i}:=\del_{v_i, \epsilon_i}$, 
the quantum partial derivative (see Definition~(\ref{qpd}))
with respect to $B_g$.

%%%%%%%%%%%%%%%%%%%%%%%%%%%%%%%%%%%%%%%%%%%%%%%%%%%%%%%
\vspace{3ex}
\begin{defn}
\label{tau}
We define a resolution converter map $\ta$ from 
the dual Koszul complex to the dual bar complex
with coefficients in $S(V)\# G$:
$$
\ta_p:\ \ 
C^p 
\rightarrow
\Hom_{\CC}(S(V)^{\ot p},S(V)\# G) \ .
$$
Let $g$ lie in $G$ with basis $B_g=\{v_1,\ldots, v_n\}$ 
of $V$ as above. Let
$
\alpha = f_g \overline{g}\otimes v_{j_1}^*\wedge\cdots\wedge v_{j_p}^*$
with 
$f_g\in S(V)$ and $1\leq j_1<\ldots<j_p \leq n$. 
Define $\ta(\alpha):S(V)^{\ot p}\rightarrow S(V)\# G$ 
by
$$
\begin{aligned}
\ta(\alpha)(& f_1\otimes \cdots\otimes f_p )
&= \Biggl(\
\prod_{k=1}^p
 {}^{s_{1}s_{2}\cdots s_{j_{k}-1}} 
(\del_{j_{k}}f_k) \Biggr) f_g\overline{g}
\ .
\end{aligned}
$$
By Theorem~\ref{psi=tau} below, $\ta$ is a cochain map. Thus
$\ta$ induces a map on the cohomology
$\HH^{\DOT}(S(V),S(V)\# G)\cong%\displaystyle
{\bigoplus_{g\in G}}
\HH^{\DOT}(S(V),S(V)\overline{g})$, 
which we denote by $\ta$ as well.
For each $g$ in $G$, let $\ta_g$ denote the restriction to 
$C^{\DOT}_g$
and the restriction to
$\HH^{\DOT}(S(V),S(V)\overline{g})$,
so that
$$
   \ta = \bigoplus_{g\in G}\ta_g.
$$
\end{defn}
%%%%%%%%%%%%%%%%%%%%%%%%%%%%%%%%%%%%%%%%%%%%%%%%%%%%%%%%%%%%
%%%%%%%%%%%%%%%%%%%%%%%%%%%
\vspace{3ex}
\begin{remark}\label{dependsonbasis}
For each $g$ in $G$, the {\em cochain} map $\ta_g=\ta_{g,B}$ 
depends on the chosen  basis $B=B_g$ of eigenvectors of $g$.  
But we shall see 
(in Corollary~\ref{IndependentOfChoice} below)
that the induced automorphism on cohomology $\cohg$
does {\em not} depend on the choice of basis.
This will imply
that as an automorphism of $\HH^{\DOT}(S(V)\# G)$, the map $\ta$
does not depend on choices of bases of $V$ used in its definition.
\end{remark}
\vspace{3ex}
%%%%%%%%%%%%%%%%%%%%%%%%%%%%%%%%%%%%%%%%%%%%%%%%%%%%%%%%%%%%%%%%
\begin{example}
Let $G=\ZZ/2\ZZ \times \ZZ/2\ZZ$ be the Klein four-group consisting of
elements $1,g,h,gh$.
Let $V=\CC^3$ with basis $v_1,v_2,v_3$ on which $G$ acts by
$ {}^gv_1=-v_1$, ${}^gv_2=v_2$, ${}^gv_3=-v_3$,
${}^hv_1=-v_1$, ${}^hv_2=-v_2$, ${}^hv_3=v_3$.
Let $\alpha=f_h\bar{h}\ot v_1^*\wedge v_2^*$ for some $f_h\in S(V)$.
Write $h=s_1s_2$, a product of reflections with 
${}^{s_1}v_1=-v_1$, ${}^{s_2}v_2=-v_2$.
Then $\ta(\alpha)$ is the function on $S(V)^{\ot 2}$ given by
$$
   \ta(\alpha) (f_1\ot f_2)=(\partial_1f_1)({}^{s_1}\partial_2f_2) f_h \bar{h}
$$
for all $f_1,f_2\in S(V)$.
For example, $\ta(\alpha)(v_1\ot v_2)=f_h\bar{h}$ while 
$\ta(\alpha)(v_2\ot v_1) =0$.
\end{example}
%%%%%%%%%%%%%%%%%%%%%%%%%%%%%%%%%%%%%%%%%%%%%%%%%%5
\vspace{3ex}
\begin{remark}
\label{zeroterm}
The map $\ta$ transforms each decomposable vector
form into a (twisted) quantum operator
characterizing the same subspace:
For the fixed basis $B_g=\{v_1,\ldots, v_n\}$ 
and $\alpha=f_g \overline{g}\otimes v_{j_1}^*\wedge\cdots\wedge v_{j_p}^*$ 
in $C^p_g$ (with $j_1<\ldots< j_p$), we have  
$$
\ta(\alpha)(v_{i_1} \otimes\cdots\ot v_{i_p})=0
\quad\quad\quad
\text{unless } i_1=j_1, \ldots, i_p=j_p\ .
$$
Generally, 
$
\ta(\alpha)(f_1\ot\cdots\ot f_p)=0
$
whenever ${\frac{\del}{\del v_{j_k}}(f_k)} =0$ for some $k$. 
\end{remark}
\vspace{3ex}

The next proposition explains how $\ta$ depends
on our choices of bases as a map on cochains.
%%%%%%%%%%%%%%%%%%%%%%%%%%%%%%%%%%%%%%%%%%%%%
\vspace{2ex}
\begin{prop}\label{changeofbasisrule}
The maps 
$\ta_{g,B}$ on cochains, for $g$ in $G$,
satisfy the following change of basis rule:
For any $a$ in $G$,
$$\,^a\ta_{g,B}= \ta_{\rule[0ex]{0ex}{1.5ex}aga^{-1},\,
  ^a\hspace{-.25ex}B}\ 
\rule[-2.5ex]{0ex}{2ex}% strut to add verticle space
 . $$
In particular, for $a$ in the centralizer $Z(g)$,
$ ^a\ta_{\rule[0ex]{0ex}{1.5ex}g,B}=\ta_{\rule[0ex]{0ex}{1.5ex}g,\, ^a\hspace{-.25ex}B}\ .$
\end{prop}
\vspace{1ex}
\begin{proof}
One may check directly from Definition~\ref{qpd} that
quantum partial differentiation obeys the following transformation law:
For all $v$ in $V$ and $\ep$ in $\CC$, 
$$
\, ^a\del_{v,\ep}=
\del_{\, ^a\hspace{-.25ex}v,\ep}\ ,
$$    
where $\, ^a\del_{v,\ep}$ differentiates with respect to a basis $B$
and $\del_{\,^av,\ep}$ with respect to $\, ^a\hspace{-.25ex}B$.

Let $B=\{v_1,\ldots, v_n\}$ be a basis of $V$ of eigenvectors of $g$
with corresponding eigenvalues $\ep_1,\ldots, \ep_n$.  
Decompose $g$ as a product of diagonal reflections
$s_i$ in $\Gl(V)$ (for $i=1,\ldots, n$) in this basis; we retain the notation before Definition~\ref{tau}.
Let $g'=aga^{-1}$, $B'=\, ^a\hspace{-.25ex}B$, 
$v_i'=\, ^av_i$, $s'_i=as_ia^{-1}$, and $\del'_i=\, ^a\del_i$.  
Then the $s'_i$ similarly decompose $g'$ in the basis $B'$ with
$\del'_i=\del_{v'_i, \ep_i}$.

Consider $\alpha=f_g\otimes v_{j_1}^*\wedge\cdots\wedge v_{j_p}^*$ in $C^p$.
For all $f_i$ in $S(V)$,
$$
\begin{aligned}
^a\bigl(\ta_{\rule[0ex]{0ex}{1.5ex}g,B}(\alpha)
\bigr)&
(f_1\otimes\ldots\otimes f_p)\\
&=
\ ^a\bigl(\ta_{\rule[0ex]{0ex}{1.5ex}g,B}(\alpha)(\, ^{a^{-1}}\scoochmore f_1\otimes\cdots\otimes
\, ^{a^{-1}}\scoochmore f_p)\bigr)\\
&=
\ ^a\biggl(\,
^{s_1\cdots s_{j_1-1}}\del_{j_1}(\, ^{a^{-1}}\scoochmore f_1)\,
\cdots\,
^{s_1\cdots s_{j_p-1}}(\del_{j_p}(\, ^{a^{-1}}\scoochmore f_p))\ f_g\ \overline{g}\biggr)\\
&=
\ 
^{s'_1\cdots s'_{j_1-1}}(^a\del_{j_1} f_1)\, 
\cdots\,
^{s'_1\cdots s'_{j_p-1}}(^a\del_{j_p}f_p)\ ^a\scoochmore f_g\ \overline{g'}\\
&=
\ 
^{s'_1\cdots s'_{j_1-1}}(\del'_{j_1}f_1)\, \cdots\,
^{s'_1\cdots s'_{j_p-1}}(\del'_{j_p}f_p)\ ^a\scoochmore f_g\ \overline{g'}\\
&=
\ta_{\rule[0ex]{0ex}{1.5ex}g',B'}(\,^a\alpha)(f_1\otimes\ldots\otimes f_p)\ ,
\end{aligned}
$$
and the result follows.
\end{proof}
%%%%%%%%%%%%%%%%%%%%%%%%%%%
\vspace{2ex}
The above proposition can also be seen  using Definition~\ref{psi} below
of the chain map
$\Psi_B$, Theorem~\ref{psi=tau} below equating $\ta_{g,B}$ and $\Psi_B$, 
and the straightforward 
fact that $\Psi_B$ has a similar change of basis property.
%%%%%%%%%%%%%%%%%%%%%%%

%%%%%%%%%%%%%%%%%%%%%%%%%%%%%%%%%%%%%%%%%%%%%%%%%%%%%%
\section{Chain maps from the bar to the Koszul resolution}
\label{chainmap}

In this section,
we define specific chain maps $\Psi_B$
from the bar resolution of $S(V)$ to its Koszul resolution
(see~(\ref{barcomplex}) and~(\ref{koszul-res2})) depending on  bases
$B$ of $V$.
By the Comparison Theorem, 
the resulting maps $(\Psi_B)^*$ on cohomology do
not depend on the choice of $B$.
In particular, we consider cohomology
with coefficients in $S(V)\overline{g}$
and write $\Psi^*_{g,B}$ for the induced map $(\Psi_{B})^*_g$ on $\cohg$.
We shall show in Theorem~\ref{psi=tau} below that $\Psi^*_{g,B} = \ta_{g,B}$ 
(recall Definition~\ref{tau}) for any choice $B$
of basis of $V$ consisting of eigenvectors of $g$ used to define both maps.
This will imply (see Corollary~\ref{IndependentOfChoice}) that as a maps on cohomology, 
$\ta_g$ and $\ta$ are automorphisms independent
of choices of bases $B$ used to define them at the cochain level. 
In this section, we work over any base field.

First, we introduce some notation.
Let $\underline{\ell}$ denote an $n$-tuple $\underline{\ell}:=(\ell_1,\ldots,\ell_n)$.
Let  $\underline{v}^{\underline{\ell}}$ be the monomial $\underline{v}^{\underline{\ell}}
:=v_1^{\ell_1}\cdots v_n^{\ell_n}$ where $v_1, \ldots, v_n$ is a chosen basis of $V$.  
Sometimes we further abbreviate a $p$-tuple $\underline{\ell}^1,\ldots, \underline{\ell}^p$
of $n$-tuples by $\underline{\ell}$ when no confusion will arise.

%%%%%%%%%%%%%%%%%%%%%
\vspace{3ex}
\begin{defn}
\label{psi}
Let $V$ be a vector space over an arbitrary field, and
let $B=\{ v_1,\ldots, v_n \}$  be a basis of $V$. 
Define an $S(V)^e$-map $\Psi_B$ from the bar resolution to the Koszul resolution, 
$
\Psi_B:S(V)^{\ot (\DOT+2)}\rightarrow S(V)^e \ot \Wedge ^{\DOT} (V)\ ,
$ 
as follows.
Let $(\Psi_B)_0$ be the identity map.
For each $p\geq 1$,   
define $(\Psi_B)_p$ by %$(\Psi_B)_k:S(V)^{\ot (k+2)}\rightarrow S(V)^e\ot \Wedge ^k(V)$ by
\begin{equation}\label{psigk}
\begin{array}{l}
(\Psi_B)_p(1\ot \underline{v}^{\underline{\ell}^1}\ot\cdots\ot \underline{v}^{\underline{\ell}^p}
   \ot 1)\\
\hspace{1.1in} =\displaystyle{\sum_{1\leq i_1<\cdots <i_p\leq n} \ 
   \sum_{0\leq a_{i_j}<\ell^j_{i_j}}
     \underline{v}^{\underline{Q}(\underline{\ell};i_1,\ldots,i_p)}\ot \underline
     {v}^{\underline{\hat{Q}}(\underline{\ell};i_1,\ldots,i_p)}\ot v_{i_1}\wedge     
    \cdots \wedge v_{i_p}},
\end{array}
\end{equation}
where the second sum ranges 
over all $a_{i_1},\ldots,a_{i_p}$ such that $0\leq a_{i_j}<\ell^j_{i_j}$
for each $j\in \{1,\ldots,p\}$
and the functions $\underline{Q}$ and $\underline{\hat{Q}}$ (indicating monomial degree)
depend also on the
choices $a_{i_j}$ (this dependence is suppressed in the notation for brevity):
$$
\underline{Q}(\underline{\ell}^1,\ldots,\underline{\ell}^p;i_1,\ldots,i_p)_i=
\left\{\begin{array}{ll}
      a_i+\ell_i^1+\cdots + \ell_i^{j-1}, & \mbox{ if } i=i_j\\
      \ell_i^1+\cdots +\ell_i^j, & \mbox{ if } i_j<i<i_{j+1},
   \end{array}\right.
$$ 
where we set $i_0=0$ and $i_{p+1}=n+1$ for convenience.
We define the $n$-tuple 
$\underline{\hat{Q}}(\underline{\ell}; i_1,\ldots, i_p)$ to be  complementary
to $\underline{Q}(\underline{\ell}; i_1,\ldots, i_p)$ in the sense that
$$
   \underline{v}^{\underline{Q}(\underline{\ell};  i_1,\ldots, i_p)}\,
 \underline{v}^{ \underline{\hat{Q}}(\underline{\ell}; i_1,\ldots, i_p)}\,  
    v_{i_1}\cdots v_{i_p} = \underline{v}^{\underline{\ell}^1}\cdots
   \underline{v}^{\underline{\ell}^p}.
$$
We simply write $\underline{\hat{Q}}$ when it is clear with respect to which
$\underline{Q}(\underline{\ell}; i_1,\ldots, i_p)$ it is complementary.
\end{defn}

\vspace{2ex}

For small values of $p$, the formula for $(\Psi_B)_p$ is less
cumbersome. In particular, for $p=1,2$, such formulas were given 
in~\cite[(4.9), (4.10)]{Witherspoon1}.
We repeat them here:
\begin{equation*}
 \Psi_1(1\ot v_1^{\ell_1}\cdots v_n^{\ell_n}\ot 1) = \sum_{i=1}^n \sum_{t=1}
   ^{\ell_i} v_i^{\ell_i-t}v_{i+1}^{\ell_{i+1}}\cdots v_n^{\ell_n}\ot 
  v_1^{\ell_1}\cdots v_{i-1}^{\ell_{i-1}} v_i^{t-1}\ot v_i,
\end{equation*}
\begin{equation*}
\Psi_2(1\ot v_1^{\ell_1}\cdots v_n^{\ell_n}\ot v_1^{m_1}\cdots v_n^{m_n}\ot 1)=
\hspace{4in}
\end{equation*}
$$\sum_{1\leq i<j\leq n} \ \sum_{r=1}^{m_j} \ \sum_{t=1}^{\ell_i}  
  v_i^{\ell_i-t} v_{i+1}^{\ell_{i+1}}\cdots v_{j-1}^{\ell_{j-1}}v_j^{\ell_j+m_j-r}v_{j+1}
  ^{\ell_{j+1}+m_{j+1}}\cdots v_n^{\ell_n+m_n}\ot 
$$

\vspace{-.15in}

$$\hspace{1.6in}v_1^{\ell_1+m_1}\cdots v_{i-1}^
 {\ell_{i-1}+m_{i-1}} v_i^{m_i+t-1} v_{i+1}^{m_{i+1}}\cdots v_{j-1}^{m_{j-1}}
  v_j^{r-1} \ot v_i\wedge v_j \ . 
$$

\vspace{3ex}
\begin{example} 
To illustrate, we compute $\Psi_2$ on a few monomials of small 
degree:
\begin{eqnarray*}
\Psi_2(1\ot v_1\ot v_2\ot 1) &=& 1\ot 1\ot v_1\wedge v_2 , \\
\Psi_2(1\ot v_1v_2\ot v_2^3\ot 1) &=& 
  (v_2^3\ot 1 + v_2^2\ot v_2 + v_2\ot v_2^2)\ot v_1\wedge v_2,\\
\Psi_2(1\ot v_1v_2\ot v_2^2v_3\ot 1) &=& (v_2^2v_3\ot 1 + v_2v_3\ot v_2)\ot
  v_1\wedge v_2 \\
  &&+ 1\ot v_1v_2^2\ot v_2\wedge v_3
   + v_2\ot v_2^2\ot v_1\wedge v_3.
\end{eqnarray*}
\end{example}
\vspace{3ex}
%%%%%%%%%%%%%%%%%%
\begin{thm}\label{psib-chain-map} For each choice of basis $B$ of $V$, the map
$\Psi_B$ of Definition~\ref{psi} is a chain map.
\end{thm}
\vspace{2ex}
%%%%%%%%%%%%%%%%%%%%%%%%%%%%%%%%%%%5
We defer the proof of Theorem~\ref{psib-chain-map} to the appendix
as it is rather technical.

%%%%%%%%%%%%%%%%%%%%%%%%%%%%%%%%%%%%%%%%%%%%%%%%%%%%%%%%%%%%%%%%%%%%%
\section{Merits of the combinatorial converter map}\label{merits}

In the previous two sections, we examined two maps $\ta_{g,B}$ and $\Psi^*_{g,B}$ 
which convert between cochain complexes:
They each transform
cochains procured from the Koszul resolution~(\ref{koszul-res2}) 
of $S(V)$ into cochains
procured from the bar resolution~(\ref{barcomplex}) of $S(V)$
(see Definitions~\ref{tau} and~\ref{psi}).
In this section, we show that the two maps
$\ta_{g,B}$ and $\Psi^*_{g,B}$ are identical on cochains, and hence 
also on cohomology,
for any $g$ in $G$ and any basis $B$
consisting of eigenvectors of $g$. 
This will imply that $\ta$ is itself a cochain map.
We deduce other salient properties of the map 
$\ta$ using this connection between $\ta$ and $\Psi$.
We take our ground field to be $\CC$ in this section.

%%%%%%%%%%%%%%%%%%%%%%%%%%%%%%%%%%5
\vspace{2ex}
\begin{theorem}
\label{psi=tau}
Let $g$ be in $G$ and let $B$ be a basis of $V$ consisting  of eigenvectors of $g$.
Then
$$\ta_{g,B}=\Psi^*_{g,B}$$ as maps on cochains.
Thus $\ta_{g,B}$ is  a cochain map.
\end{theorem}
%%%%%%%%%%%%%%%%%%%%%%%%%%%%%%%%%%%%%%
\vspace{1ex}
\begin{proof}
We check that $\ta_{g,B}$ and $\Psi_{g,B}^*$ agree on cochains:
Let $\alpha =f_g\overline{g}\ot v_{j_1}^*\wedge \cdots\wedge v_{j_p}^* $ 
be a cochain in $C^p_g$ with $f_g\in S(V)$ and $j_1<\ldots<j_p$,
where $B=\{v_1,\ldots, v_n\}$.
Let 
$f_1=\underline{v}^{\underline{\ell}^1},\ldots, f_p
=\underline{v}^{\underline{\ell}^p}$ be monomials in $S(V)$. 
Without loss of generality, it suffices to show that 
$\Psi_{g,B}^*(\alpha)$ and $\ta_{g,B}(\alpha)$ 
agree on $f_1\ot \cdots\ot f_p$, since such elements form 
a basis for $S(V)^{\ot p}$.  By Definition~\ref{psi},
$$
\begin{aligned}
\Psi_{g,B}^*(\alpha)(&f_1\ot    \cdots\ot f_p)\\
&=\ \alpha(\Psi_{g,B}(f_1\ot\cdots\ot f_p))\\
&=\ \alpha(\Psi_{g,B}( v_{1}^{\ell^1_1}\cdots v_{n}^{\ell^1_n}
   \ot\cdots\ot  v_{1}^{\ell^p_1}\cdots v_{n}^{\ell^p_n}))\\
&=\ \alpha\biggl(\
\sum_{1\leq i_1<\cdots <i_p\leq n} \ 
   \sum_{0\,\leq\, a_{i_k}<\,\ell^k_{i_k}}
     \underline{v}^{\underline{Q}(\underline{\ell};i_1,\ldots,i_p)}\ot \underline
     {v}^{\underline{\hat{Q}}(\underline{\ell};i_1,\ldots,i_p)}\ot 
      v_{i_1}\wedge \cdots \wedge v_{i_p}\biggr)\ .
\end{aligned}
$$
Since  $\alpha$ has exterior part
$v_{j_1}^*\wedge \cdots\wedge v_{j_p}^*$, 
each summand is zero save one (the summand with $i_k=j_k$ for $k=1,\ldots p$).
Then

\vspace{1ex}
$$
\begin{aligned}
&\Psi^*_{g,B}(\alpha)(f_1\ot\cdots\ot f_p)\\
&=
\sum_{0\,\leq\, a_{j_k}<\,\ell^k_{j_k}}
     \underline{v}^{\underline{Q}(\underline{\ell};j_1,\ldots,j_p)}
    \ \ f_g \overline{g}\ \  
     \underline{v}^{\underline{\hat{Q}}(\underline{\ell};j_1,\ldots,j_p)}\\
&=
 \sum_{0\,\leq\, a_{j_k}<\,\ell^k_{j_k}}
   \biggl(\ \prod_{t=1}^p \epsilon_{j_t}^{(\ell_{j_t}^{t}-a_{j_t}-1)+\ell_{j_t}^{t+1}
    \cdots + \ell_{j_t}^p} \prod_{i_{t-1} < i <i_t} \epsilon_i^{\ell_i^t +
    \cdots + \ell_i^p}
      \biggr)\ 
     \underline{v}^{\underline{Q}(\underline{\ell};j_1,\ldots,j_p)}\underline
     {v}^{\underline{\hat{Q}}(\underline{\ell};j_1,\ldots,j_p)}
         \ \   f_g\overline{g}.
\end{aligned}
$$

\vspace{2ex}

\noindent
Recall that by definition,  
$\underline{v}^{\underline{Q}(\underline{\ell};j_1,\ldots,j_p)}\underline
     {v}^{\underline{\hat{Q}}(\underline{\ell};j_1,\ldots,j_p)} v_{j_1}
  \cdots v_{j_p} = \underline{v}^{\underline{\ell}^1}\cdots 
  \underline{v}^{\underline{\ell}^p}$.
Therefore the factor 
$\underline{v}^{\underline{Q}(\underline{\ell};j_1,\ldots,j_p)}\underline
     {v}^{\underline{\hat{Q}}(\underline{\ell};j_1,\ldots,j_p)} f_g\overline{g}$
in each term of the above sum does not depend on the values of $a_{j_k}$,
and thus we may move the summation symbol inside the parentheses.
Simplifying, we obtain

\vspace{1ex}
$$
\begin{aligned}
&
 \biggl(\ \prod_{t=1}^p  [\ell_{j_t}^t]_{\epsilon_{j_t}} 
       \epsilon_{j_t}^{\ell_{j_t}^{t+1}+\cdots +\ell^p_{j_t}}
    \prod_{i_{t-1}<i<i_t}\epsilon_i^{\ell_i^t+\cdots +\ell_i^p}
    \biggr)\ \
  \underline{v}^{\underline{Q}(\underline{\ell};j_1,\ldots,j_p)}\underline
     {v}^{\underline{\hat{Q}}(\underline{\ell};j_1,\ldots,j_p)}
      f_g\overline{g}\\ 
&= \Bigl(\
  \prod_{k=1}^p  {}^{s_{1}\cdots s_{j_{k}-1}} (\partial_{j_k}f_k)
  \Bigr) \ \ f_g\overline{g} \ \\
&=
\ta_{\rule[0ex]{0ex}{1.5ex}g,B}(\alpha)(     f_1\ot\cdots\ot f_p     ),
\end{aligned}
$$

\vspace{2ex}

\noindent
by Definition~\ref{tau}.
Hence $\ta_{g,B}=\Psi_{g,B}^*$ as maps on cochains.
\end{proof}
\vspace{2ex}

As a consequence, we obtain a statement about any chain map $\Psi$
from the bar resolution to the Koszul resolution, at the level of cohomology,
which has further implications for the map $\ta$:

%%%%%%%%%%%%%%%%%%%%%%%%%%%%%%%%%%%
\vspace{2ex}
\begin{cor}
\label{tau=anypsi}
Let $\Psi$ be any chain map from the bar resolution~(\ref{barcomplex})
to the Koszul resolution~(\ref{koszul-res2}) for $S(V)$. Then:
\begin{enumerate}
\item
$\ta_g=\Psi_g^*$ as maps on the cohomology
$\HH^{\DOT}(S(V), S(V)\bar{g})$, for all $g$ in $G$.
\item
$\ta=\Psi^*$
as maps on 
$$ 
  \HH^{\DOT}(S(V),S(V)\# G)\cong \bigoplus_{g\in G}\, \cohg,
$$ 
and on its  $G$-invariant subalgebra, 
$$ 
  \HH^{\DOT}(S(V),S(V)\# G)^G\cong  \HH^{\DOT}(S(V)\#G).
$$
\end{enumerate}
\end{cor}
\vspace{1ex}

\begin{proof}
We constructed a specific choice of chain map $\Psi_B$  in
Definition~\ref{psi} above
from the bar to the Koszul resolution of $S(V)$.
Since $\Psi$ and $\Psi_B$ are homotopic
by the Comparison Theorem, $\Psi^*_g=\Psi_{g,B}^*$ as maps on cohomology
$\cohg$.  But $\Psi_{g,B}^*=\ta_{g,B}$ 
for any choice of $g$ and $B$ by Theorem~\ref{psi=tau}, and 
hence $\Psi^*=\ta$.
By Proposition~\ref{psiInvariant}, these maps preserve
$G$-invariant subspaces, and so $\Psi^*=\ta$ on
$\HHD(S(V)\# G)$ as well.
\end{proof}

%%%%%%%%%%%%%%%%%%%%%%%%%%%%%%%%%%%%%%%%%%%%%%%%%%%%%%5
\vspace{2ex}
\begin{cor}
\label{IndependentOfChoice}
Let $g\in G$.
On the cohomology $\HH^{\DOT}(S(V), S(V)\bar{g})$, 
the map $\ta_g=\ta_{g,B}$ is 
independent of choice of basis $B$ of eigenvectors of $g$ used in its definition.
Hence, as a map on the cohomologies $\HH^{\DOT}(S(V),S(V)\# G)$ and $\HH^{\DOT}(S(V)\# G)$, 
$\ta$ is independent of the choices of bases used in its definition.
\end{cor}
%%%%%%%%%%%%%%%%%%%%%%%%%%%%%%%%%%%%%%%%%%%%%%%%%%%%
\vspace{1ex}
\begin{proof}
By Corollary~\ref{tau=anypsi}, $\ta_g=\Psi_g^*$ on cohomology 
for any choice of chain map $\Psi$ from the bar complex to the Koszul complex
of $S(V)$, 
independent of the choice of 
basis of eigenvectors of $g$ used to define $\ta_g$.  Hence,
$\ta$ is independent of choices of bases. 
\end{proof}
%%%%%%%%%%%%%%%%%%%%%%%%%%%%%%%%5
\vspace{2ex}
%%%%%%%%%%%%%%%%%%%%%%%%%%%%%%%5
\begin{cor}\label{inverseisos}
The maps $\ta$ and $\Phi^*$ are inverse isomorphisms 
on the cohomology
$\HH^{\DOT}(S(V),S(V)\# G)$ and on its $G$-invariant subalgebra
$\HH^{\DOT}(S(V),S(V)\# G)^G\cong \HH^{\DOT}(S(V)\# G)$. 
\end{cor}
%%%%%%%%%%%%%%%%%%%%%%%%%%%%%%%%%%%
\vspace{1ex}
\begin{proof}
Again, in Theorem~\ref{tau=anypsi}, we found that $\ta=\Psi^*$ on cohomology
for any chain map $\Psi$ from the bar to the Koszul complex.  
But any such $\Psi$ induces an automorphism on cohomology inverse
to $\Phi^*$ by the Comparison Theorem.
\end{proof}
%%%%%%%%%%%%%%%%%%%%%%%%%%%%%%%%%%%%%%%%%%%%%%%%%%%%%%%
\vspace{2ex}
The above corollary actually follows from a stronger fact:
$\ta$ is a right-sided inverse to $\Phi^*$ on cochains, not
merely on cohomology, for any choice of bases $\{B_g\}_{g\in G}$
defining $\ta$.
Indeed, a calculation shows directly 
that $\Phi^*\ta=1$ on cochains $C^{\DOT}$.
We can see this fact yet another way.
One can check that $\Psi_B\Phi=1$ on chains, 
and therefore $\Phi_g^*(\Psi_B)_g^*=1$
on cochains, for every $B$ and $g$.  
In Theorem~\ref{psi=tau}, we saw
that $\ta_{g,B}=(\Psi_B)_g^*=\Psi_{g,B}^*$ as maps on cochains,
for all $g$ in $G$ and for any basis $B$ of eigenvectors of $g$, 
and hence $\Phi^* \ta=1$ as a map on cochains.
%%%%%%%%%%%%%%%%%%%%%%%%%%%%%%%%%%%%%%%%%%%%%%%

%%%%%%%%%%%%%%%%%%%%%%%%%%%%%%%%%%%%%%%%%%%%%%%
%%%%%%%%%%%%%%%%%%%%%%%%%%%%%%%%%%%%%%%%%%%%%%%
%%%%%%%%%%%%%%%%%%%%%%%%%%%%%%%%%%%%%%%%%%%%%%%
\section{Hochschild homology}\label{intheliterature}
Our chain maps $\Psi_B$ of Definition~\ref{psi}
are useful in settings other than the cohomology of $S(V)\# G$.
In this section, we obtain induced maps on 
Hochschild homology, and compare our induced maps on
homology and cohomology with those in the literature.
The Hochschild-Kostant-Rosenberg Theorem states that
for smooth commutative algebras, Hochschild homology is isomorphic to 
the module of differential forms (i.e., the exterior algebra
generated by the K\"ahler differentials);
e.g., see~\cite[\S9.4.2]{Weibel}.
For noncommutative algebras, Hochschild homology provides a generalization 
of the notion of ``differential forms''.
It is interesting to note that for some types of algebras
(in particular for $S(V)\# G$), Hochschild homology and cohomology are
dual
(see van den Bergh~\cite{vandenBergh} for the general theory
and Farinati~\cite{Farinati} for the case $S(V)\# G$).
In this section, we work over an arbitrary field initially,
then over $\CC$ in Theorem~\ref{induced-homology}.

Let $M$ be any $S(V)^e$-module.
Then 
$\Psi_B$ induces an isomorphism on Hochschild homology $\HH_{\DOT}(S(V),M)
:= \Tor^{S(V)^e}_{\DOT}(S(V), M)$
and on Hochschild cohomology $\HH^{\DOT}(S(V), M):=\Ext_{S(V)^e}^{\DOT}(S(V),M)$ 
by applying the functors 
$M\ot _{S(V)^e} - $ and $\Hom_{S(V)^e}( - , M)$, respectively, to the 
bar resolution~(\ref{barcomplex}) and to the Koszul resolution~(\ref{koszul-res2}).
This approach to obtaining maps on homology and cohomology has advantages over
previous approaches in the literature which we explain now.

We obtain a map $(\Psi_B)_*$ on Hochschild homology $\HH_{\DOT}(S(V)):=
\HH_{\DOT}(S(V), S(V))$ by setting $M=S(V)$.
At the chain level, 
$$(\Psi_B)_*:S(V)\ot S(V)^{\ot \DOT}\longrightarrow S(V)\ot \Wedge^{\DOT}(V)\ .$$
A computation similar
to that in the proof of Theorem~\ref{psi=tau} yields the following
explicit formula for $(\Psi_B)_*$, 
valid over any ground field:
%%%%%%%%%%%%%%%%%%%%%%%%%%%%%%%%%%%%%%%%%%%%%%%%%%
\vspace{2ex}
\begin{thm}\label{Halbout-map}
Let $B=\{v_1,\ldots,v_n\}$ be a basis of $V$.
Then as an automorphism on $\HH_{\DOT}(S(V))$ at the chain level,
\begin{equation}
    (\Psi_B)_* ( f_0\ot f_1\ot\cdots\ot f_p) =
   \sum_{1\leq i_1 < \cdots < i_p\leq n} f_0\ \frac{\partial f_1}{\partial v_{i_1}}
   \cdots \frac{\partial f_p}{\partial v_{i_p}} \ot v_{i_1}\wedge\cdots\wedge v_{i_p}
\end{equation}
for all $f_0,f_1,\ldots,f_p\in S(V)$.
\end{thm}
%%%%%%%%%%%%%%%%%%%%%%%%%%%%%%%%%%%%%%%%%%%%%%%%
\vspace{1ex}
\begin{proof}
Without loss of generality assume that $f_1,\ldots,f_p$ are monomials, say
each $f_k=\underline{v}^{\underline{\ell}^k}
=v_1^{\ell_1^k}\cdots v_n^{\ell_n^k}$ for some $n$-tuple $\underline{\ell}^k=(\ell_1^k,\ldots,\ell_n^k)$.
Then
\begin{eqnarray*}
 && \hspace{-1.6cm} (\Psi_B)_*(f_0\ot \underline{v}^{\underline{\ell}^1}\ot\cdots\ot 
        \underline{v}^{\underline{\ell}^p})\\
     & = & f_0 \Psi_B(1\ot 
     \underline{v}^{\underline{\ell}^1}\ot\cdots\ot \underline{v}^{\underline{\ell}^p})\\
    &=& f_0 \left(\sum_{1\leq i_1<\cdots<i_p\leq n} \ 
     \sum_{0\leq a_{i_j}<\ell^j_{i_j}} \underline{v}^{\underline{Q}(\underline{\ell};
    i_1,\ldots,i_p)}\ot \underline{v}^{\underline{\hat{Q}}(\underline{\ell};
     i_1,\ldots,i_p)}\ot v_{i_1}\wedge\cdots\wedge v_{i_p}\right)\\
   &= &  \sum_{1\leq i_1<\cdots<i_p\leq n} \ 
     \sum_{0\leq a_{i_j}<\ell^j_{i_j}} \underline{v}^{\underline{Q}(\underline{\ell};
    i_1,\ldots,i_p)}\ f_0 \ \underline{v}^{\underline{\hat{Q}}(\underline{\ell};
     i_1,\ldots,i_p)}\ot v_{i_1}\wedge\cdots\wedge v_{i_p}\\
&= &  \sum_{1\leq i_1<\cdots<i_p\leq n} \ 
     \sum_{0\leq a_{i_j}<\ell^j_{i_j}} f_0\underline{v}^{\underline{Q}(\underline{\ell};
    i_1,\ldots,i_p)}\underline{v}^{\underline{\hat{Q}}(\underline{\ell};
     i_1,\ldots,i_p)}\ot v_{i_1}\wedge\cdots\wedge v_{i_p}\\
&= &  \sum_{1\leq i_1<\cdots<i_p\leq n} \ 
     \sum_{0\leq a_{i_j}<\ell^j_{i_j}} f_0 \underline{v}^{\underline{\ell}^1}\cdots
     \underline{v}^{\underline{\ell}^p} v_{i_1}^{-1}\cdots v_{i_p}^{-1}
    \ot v_{i_1}\wedge\cdots\wedge v_{i_p}\\
&= &  \sum_{1\leq i_1<\cdots<i_p\leq n} \ 
     \ell_{i_1}^1\cdots \ell_{i_p}^p f_0 \underline{v}^{\underline{\ell}^1}\cdots
     \underline{v}^{\underline{\ell}^p} v_{i_1}^{-1}\cdots v_{i_p}^{-1}
     \ot v_{i_1}\wedge\cdots\wedge v_{i_p},
\end{eqnarray*}
where the product  $\underline{v}^{\underline{\ell}^1}\cdots
     \underline{v}^{\underline{\ell}^p} v_{i_1}^{-1}\cdots v_{i_p}^{-1}$ is
computed in the ring of Laurent polynomials in $v_1,\ldots,v_n$.
(Since $0\leq a_{i_j}<\ell_{i_j}^j$, the result lies in $S(V)$ when the 
corresponding sum is nonempty.)
The expression above is precisely that claimed in the theorem.
\end{proof}
\vspace{2ex}

In case the ground field is $\CC$ or ${\mathbb R}$, by the above
theorem, our map $(\Psi_B)_*$ is precisely
the map $J$ of Halbout~\cite{Halbout}.
Halbout gave an explicit homotopy $s$ showing that $J$ is a homotopy inverse to the
canonical embedding of  the
de Rham complex into the Hochschild complex.
In contrast, we see immediately that $(\Psi_B)_*$ induces 
an isomorphism on homology since $\Psi_B$ is itself a chain map.

For comparison, we give the map on Hochschild cohomology $\HHD(S(V))$;
this is simply the case $g=1$ of Definition~\ref{tau}, by Theorem~\ref{psi=tau}:

\vspace{2ex}
%%%%%%%%%%%%%%%%%%%%%%%%%%%%%
\begin{thm}
Let $B=\{v_1,\ldots,v_n\}$ be a basis of $V$.
Then as an automorphism on $\HHD(S(V))$ at the chain level,
$$
  (\Psi_B)^*(\alpha)(f_1\ot\cdots\ot f_p) = f_0\
   \frac{\partial f_1}{\partial v_{j_1}}\cdots \frac{\partial f_p}{\partial v_{j_p}}
$$
when $\alpha = f_0\ot v_{j_1}^*\wedge\cdots\wedge v_{j_p}^*\in \Hom_{\CC}(S(V)\ot
\Wedge^pV , S(V))$,
$f_1,\ldots,f_p\in S(V)$.
\end{thm}
\vspace{2ex}
%%%%%%%%%%%%%%%%%%%%%%%%%%%%%

Now we restrict our choice of field again to $\CC$.
Let $M=S(V)\# G$, and note that Hochschild homology decomposes
just as does Hochschild cohomology:
$$
  \HH_{\DOT}(S(V)\# G)\cong \HH_{\DOT}(S(V),S(V)\# G)^G\cong
   \left( \bigoplus_{g\in G} \HH_{\DOT}(S(V),S(V)\bar{g})\right)^{G}
$$
(see~\cite{Farinati,Stefan}).
Thus one is interested in the components $\HH_{\DOT}(S(V),S(V)\bar{g})
=\Tor_{\DOT}^{S(V)^e} (S(V),S(V)\bar{g})$, for each $g$ in $G$.
A calculation similar to that in the proof of Theorem~\ref{psi=tau}
yields the explicit formula in the next theorem for the induced map
$$(\Psi_B)_*:S(V)\overline{g}\ot S(V)^{\ot \DOT}
\rightarrow
S(V)\bar{g}\ot \Wedge^{\DOT}(V)\ .
$$
Note that quantum differential operators surface (compare with 
Definition~\ref{tau} of $\ta$, which 
is equal to $\Psi_B^*$ by Theorem~\ref{psi=tau}). 
We have not found these maps in the literature on Hochschild homology.

For $g$ in $G$, let $B=\{v_1,\ldots,v_n\}$ be a basis of $V$ consisting
of eigenvectors of $g$ with corresponding eigenvalues $\ep_1,\ldots, \ep_n$.
Write $g=s_1\cdots s_n$ where 
${}^{s_i}v_j=v_j$ for $j\neq i$ and $ {}^{s_i} v_i = \epsilon_i v_i$.
Recall the quantum partial differential operators $\partial_i:=\partial_{v_i,
  \ep_i}$ of 
Definition~(\ref{qpd}).

%%%%%%%%%%%%%%%%%%%%%%%%%%%%%%%
\vspace{2ex}
\begin{thm}\label{induced-homology}
Let $g\in G$ and let $B=\{v_1,\ldots,v_n\}$ be a basis of $V$ consisting of eigenvectors for $g$.
Then as an automorphism on $\HH_{\DOT}(S(V),S(V)\bar{g})$ at the chain level, 
$$
  (\Psi_B)_*(f_0\bar{g}\ot f_1\ot\cdots\ot f_p)=
   \sum_{1\leq i_1<\cdots <i_p \leq n} f_0 \left(\prod_{k=1}^p
      {}^{s_1s_2\cdots s_{i_{k-1}}} (\partial_{i_k} f_k)\right) \bar{g} \ot
    v_{i_1}\wedge \cdots \wedge v_{i_p}
$$
for all $f_0,f_1,\ldots,f_p\in S(V)$.
\end{thm} 
\vspace{2ex}
%%%%%%%%%%%%%%%%%%%%%%%%%%%%%

We make a few final comments about the appearance of our chain maps $\Psi_B$
in Hochschild cohomology.
Again let $M=S(V)\# G$, and consider the map $(\Psi_B)^*$ 
on the Hochschild cohomology $\HH^{\DOT}(S(V),S(V)\# G)$
for any basis $B$ of $V$.
We observed
(as a consequence of Theorem~\ref{psi=tau} and Definition~\ref{tau})
that $(\Psi_B)^*$ is given by quantum partial differential operators.
The reader should compare with maps given in 
Halbout and Tang~\cite{HalboutTang}:
They define functions directly on cochain complexes 
(without first defining chain maps on resolutions)
and then must prove that these functions are cochain maps.
Again, our approach presents an advantage:
We instead define one primitive chain map $\Psi_B$
from which induced cochain maps effortlessly spring.
For example, 
$(\Psi_B)^* = \ta$ is automatically a cochain map
since $\Psi_B$ is a chain map by Theorem~\ref{psib-chain-map}. 
The reader is cautioned that Halbout and Tang (in~\cite{HalboutTang})
work only over $\mathbb{R}$, in which case Hochschild cohomology
has a specialized description ($V$ and $V^*$ are $G$-isomorphic in that setting, simplifying
some aspects of homology and cohomology).

%%%%%%%%%%%%%%%%%%%%%%%%%%%%%%%%%%%%%%%%%%%%%%%%%%%%%%%%%%%%%%
%%%%%%%%%%%%%%%%%%%%%%%%%%%% SECTION %%%%%%%%%%%%%%%%%%%%%%%%%
%%%%%%%%%%%%%%%%%%%%%%%%%%%%%%%%%%%%%%%%%%%%%%%%%%%%%%%%%%%%%%
\section*{Appendix: Proof of Theorem~\ref{psib-chain-map}}\label{sec:psigk}

\maketitle

Let $V$ be a finite dimensional vector space over any field.
Fix a basis $B=\{v_1,\ldots,v_n\}$ of $V$.
Recall Definition~\ref{psi} of the linear map $\Psi=\Psi_B$ from the
bar resolution~(\ref{barcomplex}) to the Koszul resolution~(\ref{koszul-res2}) of $S(V)$.
We now prove that $\Psi$ is a chain map:
We show that
$\Psi_{p-1}\delta_p = d_p \Psi_p$ for all $p\geq 1$.

A straightforward but tedious calculation shows that
$\Psi_0\delta_1 = d_1\Psi_1$, and 
we assume from now on that $p\geq 2$.

We first compute $d_p\Psi_p(1\ot\underline{v}^{\underline{\ell}^1}\ot\cdots
\ot \underline{v}^{\underline{\ell}^p}\ot 1)$.
For each $j\in\NN$,
let $\delta[j] : \NN\rightarrow \{0,1\}$
be the Kronecker delta function defined by
$$
  (\delta[j])_i = \delta[j](i) = \left\{\begin{array}{ll}
     1, & \mbox{ if } i=j\\  0, & \mbox{ if } i\neq j.\end{array}\right.
$$
Then
\begin{eqnarray*}
 &&\hspace{-1cm} d_p\Psi_p(1\ot\underline{v}^{\underline{\ell}^1}\ot\cdots
  \ot \underline{v}^{\underline{\ell}^p}\ot 1)\hspace{1in}\\
 &=& d_p(\displaystyle{\sum_{1\leq i_1<\cdots <i_p\leq n} \ 
   \sum_{0\,\leq\, a_{i_j}<\,\ell^j_{i_j}}
     \underline{v}^{\underline{Q}(\underline{\ell};i_1,\ldots,i_p)}\ot \underline
     {v}^{\underline{\hat{Q}}(\underline{\ell};i_1,\ldots,i_p)}\ot v_{i_1}\wedge     
    \cdots \wedge v_{i_p}})\\
  &=& \sum_{1\leq i_1<\cdots <i_p\leq n} \ 
   \sum_{0\,\leq\, a_{i_j}<\,\ell^j_{i_j}}
    \ \sum_{m=1}^p (-1)^{m+1}\left(\underline{v}^{\underline{Q}
       +\delta[i_m]}\ot\underline{v}^{\underline{\hat{Q}}}\ot v_{i_1}\wedge\cdots\wedge
          \hat{v}_{i_m}\wedge\cdots\wedge v_{i_p}\right. \\
   &&   \hspace{5.5cm}- \ \left.\underline{v}^{\underline{Q}
       }\ot\underline{v}^{\underline{\hat{Q}}+\delta[i_m]}\ot v_{i_1}\wedge\cdots\wedge
          \hat{v}_{i_m}\wedge\cdots\wedge v_{i_p}\right),
\end{eqnarray*}
where $\underline{Q}=\underline{Q}(\underline{\ell}; i_1,\ldots,i_p)$
and ${\underline{\hat{Q}}}={\underline{\hat{Q}}}(\underline{\ell};i_1,\ldots,i_p)$
are determined by the $a_{i_j}$ as in Definition~\ref{psi}.
Now fix $m$ in the above expression. The factors
$\underline{v}^{\underline{Q}+\delta[i_m]}$
and $\underline{v}^{\underline{Q}}$ differ only in the
power of $v_{i_m}$.
In the sum, the power $a_{i_m}$ ranges over the set $\{0,\ldots,\ell^m_{i_m}-1\}$, and
thus the corresponding terms 
cancel except for the first term when $a_{i_m}=\ell^m_{i_m}-1$ and
the second term when $a_{i_m}=0$. After all such cancellations, for each
$m=1,\ldots,p$, what remains is
$$
\begin{aligned}
&  \sum_{1\leq i_1<\cdots <i_p\leq n} \ 
   \sum_{m=1}^p  \
   \sum_{0\,\leq\, a_{i_j}<\,\ell^j_{i_j}}
    \! (-1)^{m+1}(\underline{v}^{\underline{R}^m(\underline{\ell};i_1,\ldots,i_p)}\ot
     \underline{v}^{\underline{\hat{R}}^m}\ot v_{i_1}\wedge\cdots\wedge \hat{v}_{i_m}
     \wedge\cdots\wedge v_{i_p} \\
  &  \hspace{2.25in}- \ \underline{v}^{\underline{S}^m(\underline{\ell};i_1,\ldots,i_p)}\ot
     \underline{v}^{\underline{\hat{S}}^m}\ot v_{i_1}\wedge\cdots\wedge \hat{v}_{i_m}
     \wedge\cdots\wedge v_{i_p}),
\end{aligned}
$$
where the rightmost sum is over $a_i,\ldots, \hat{a}_{i_m},\ldots,a_{i_p}$ and
(abusing notation, as 
 $\underline{R},\underline{S}$ depend on fewer $a_{i_j}$'s than
$\underline{Q}$)
\begin{eqnarray*}
  \underline{R}^m(\underline{\ell};i_1,\ldots,i_p)_i &=& 
   \left\{\begin{array}{ll} \ell^1_{i_m}+\cdots + \ell^m_{i_m}, & \mbox{ if }i=i_m\\ 
           \underline{Q}(\underline{\ell};i_1,\ldots,i_p)_i, & \mbox{ if } i\neq i_m
          \end{array}\right. \\
\mbox{and } \ \underline{S}^m(\underline{\ell};i_1,\ldots,i_p)_i &=& 
   \left\{\begin{array}{ll} \ell^1_{i_m}+\cdots + \ell^{m-1}_{i_m}, & \mbox{ if }i=i_m\\ 
           \underline{Q}(\underline{\ell};i_1,\ldots,i_p)_i, & \mbox{ if } i\neq i_m
          \end{array}\right.,
\end{eqnarray*}
and $\underline{\hat{R}}^m, \underline{\hat{S}}^m$ are defined by the equations
\begin{eqnarray*}
  \underline{v}^{ \underline{R}^m(\underline{\ell};i_1,\ldots,i_p)}
  \underline{v}^{\underline{\hat{R}}^m
  (\underline{\ell};i_1,\ldots,i_p)} v_{i_1}\cdots\hat{v}_{i_m}\cdots v_{i_p} &=&
    \underline{v}^{\underline{\ell}^1}\cdots\underline{v}^{\underline{\ell}^p},\\
  \underline{v}^{\underline{S}^m(\underline{\ell};i_1,\ldots,i_p)}
  \underline{v}^{\underline{\hat{S}}^m
  (\underline{\ell};i_1,\ldots,i_p)} v_{i_1}\cdots\hat{v}_{i_m}\cdots v_{i_p} &=&
    \underline{v}^{\underline{\ell}^1}\cdots\underline{v}^{\underline{\ell}^p}.
\end{eqnarray*}

Consider the leftmost sum over $1\leq i_1<\cdots <i_p\leq n$. 
If we replace a given $i_m$ in $\underline{S}^m$ by $i_m+1$ in $\underline{R}^m$
(provided $i_m+1 <i_{m+1}$), keeping the others fixed,
then
$\underline{S}^m(\underline{\ell};i_1,\ldots,i_p)=
\underline{R}^m(\underline{\ell};i_1,\ldots,i_{m-1},i_m+1,i_{m+1},\ldots,i_p)$.
We thus have further cancellation, with the remaining terms 
coming from the first summand when $i_m=i_{m-1}+1$ and the second summand
when $i_m=i_{m+1}-1$:
\begin{eqnarray*}
\hspace{-1cm}&&  \sum_{m=1}^p \ \
\sum_{1\leq i_1<\cdots<\hat{i}_m<\cdots <i_p\leq n} \ \  
   \sum_{0\,\leq\, a_{i_j}<\,\ell^j_{i_j}}
    \! \! (-1)^{m+1}\\
&&\hspace{1.7cm} (\underline{v}^{\underline{R}^m(\underline{\ell};i_1,\ldots,i_{m-1},i_{m-1}+1,i_{m+1},
\ldots,i_p)}\! \ot\! \underline{v}^{\underline{\hat{R}}^m}\! \ot \!
   v_{i_1}\wedge\cdots\wedge\hat{v}_{i_m}\wedge
     \cdots\wedge v_{i_p}\\
&& \hspace{1.5cm}
  - \ \underline{v}^{\underline{S}^m(\underline{\ell};i_1,\ldots,i_{m-1},i_{m+1}-1,i_{m+1},
\ldots,i_p)}\!\ot\!\underline{v}^{\underline{\hat{S}}^m}\!\ot\! v_{i_1}\wedge
  \cdots\wedge\hat{v}_{i_m}\wedge
\cdots\wedge v_{i_p}),
\end{eqnarray*}
where the rightmost sum is over all such $a_{i_1},\ldots, \hat{a}_{i_m},\ldots, a_{i_p}$.

Now consider the middle sum above ranging 
over all $1\leq i_1<\cdots<\hat{i}_m<\cdots<i_p\leq n$.
If $m=1$, this sum does not include $i_2=1$, due to the left out entry
$\hat{i}_m$. Similarly, if $m=p$, this sum does not include $i_{p-1}=n$.
For the sake of later comparison, we add and subtract terms in the $m=1$
summand, corresponding to $i_2=1$, and in the
$m=p$ summand, corresponding to $i_{p-1}=n$. 
These added and subtracted terms may be written with either notation,
 $\underline{R}$ or $\underline{S}$, so that the result looks the same as above
except that now we include summands corresponding to $m=1,i_2=1$ and to $m=p,i_{p-1}=n$.

We next combine some of the terms. 
Consider the terms arising from a pair of subsequent indices $m$ and $m+1$
in the leftmost sum.  We pair each summand of type $\underline{S}^m$
with a summand of type $\underline{R}^{m+1}$.
Fix an integer $i$
and collect those summands
(in the $m$-th sum) with 
$\underline{S}^m$-exponent for which $i_{m+1} = i$ and those
summands (in the $(m+1)$-st sum) with $\underline{R}^{m+1}$-exponent for which $i_m =i$.  
Note that all these summands share the same sign.
We compare the exponents
$\underline{S}^m(\underline{\ell};
 i_1,\ldots,i_{m-1},i_{m+1}-1,i_{m+1},\ldots,i_p)$
and
$\underline{R}^{m+1}(\underline{\ell};i_1,\ldots,i_{m},i_{m}+1,i_{m+2},
\ldots,i_p)$ when $i_{m}=i=i_{m+1}$ and see that 
the power of $v_i$ ranges
from $\ell_i^1+\cdots +\ell_i^{m-1}$ to $\ell_i^1+\cdots +\ell_i^{m-1} +\ell_{i}^m-1$ 
and then again from $\ell_i^1+\cdots +\ell_i^{m-1}+\ell_{i}^m$
to $\ell_i^1+\cdots +\ell_{i}^m+\ell_{i}^{m+1}-1$ in this collection.
Hence, we can simply rewrite the partial sum over this collection
using the exponent
$\underline{Q}(\underline{\ell}^1,\ldots,\underline{\ell}^{m-1},
\underline{\ell}^m+\underline{\ell}^{m+1},\underline{\ell}^{m+2},\ldots,
\underline{\ell}^p; i_1,\ldots,\hat{i}_{m+1},\ldots,i_p)$ instead.
We obtain the following, in which the $m=1$ (unmatched $\underline{R}^1$)
and $m=p$ (unmatched $\underline{S}^{p}$) sums have
been singled out:
$$
\begin{aligned}
&\sum_{1\leq i_2<\cdots <i_p\leq n} \  \
   \sum_{\substack{0\,\leq\, a_{i_j}<\, \ell^j_{i_j}\\(
{\mbox{\scriptsize{for} }}
j\in \{2,\ldots,p\})}}
  \underline{v}^{\underline{Q}(\underline{\ell}^2,\cdots,\underline{\ell}^p;
   i_2,\cdots,i_p)+\underline{\ell}^1}\ot \underline{v}^{\underline{\hat{Q}}}
  \ot v_{i_2}\wedge\cdots\wedge v_{i_p}\\
&+\sum_{m=1}^{p-1}(-1)^{m}
   \sum_{1\leq i_1<\cdots<\hat{i}_m<\cdots <i_p\leq n} \  
   \sum_{\substack{ 0\, \leq\, a_{i_j}<\, \ell^{j}_{i_j} \ (
{\mbox{\scriptsize{for} }}
j\in \{1,\ldots, m-1\})\\
\rule[-1.5ex]{0ex}{3.5ex} %strut
0\,\leq\, a_{i_{m+1}}\leq\, \ell_{i_{m+1}}^m+\ell_{i_{m+1}}^{m+1}-1\\
     0\,\leq\, a_{i_j}<\,\ell^{j+1}_{i_j} \ (
{\mbox{\scriptsize{for} }}
j\in \{m+2,\ldots,p\})}}\\
&\hspace{.5in}
\rule[-2ex]{0ex}{6ex} %strut
  \underline{v}^{\underline{Q}(\underline{\ell}^1,\ldots,\underline{\ell}^{m-1},
\underline{\ell}^m+\underline{\ell}^{m+1},\underline{\ell}^{m+2},\ldots,
\underline{\ell}^p; i_1,\ldots,\hat{i}_{m},\ldots,i_p)}\ot
\underline{v}^{\underline{\hat{Q}}}\ot v_{i_1}\wedge\cdots\wedge\hat{v}_m\wedge
\cdots\wedge v_{i_p}\\
&+(-1)^{p} \!  
\sum_{1\leq i_1<\cdots <i_{p-1}\leq n} 
   \sum_{\substack{0\,\leq\, a_{i_j}<\,\ell^j_{i_j}\\ (
{\mbox{\scriptsize{for} }}
j\in \{1,\ldots,p-1\})}}
  \underline{v}^{\underline{Q}(\underline{\ell}^1,\cdots,\underline{\ell}^{p-1};
   i_2,\cdots,i_{p-1})}\ot \underline{v}^{\underline{\hat{Q}} +\underline{\ell}^p}
  \ot v_{i_1}\wedge\cdots\wedge v_{i_{p-1}}.
\end{aligned}
$$

Now relabel indices
so that each sum is taken over
$1\leq i_1<\cdots<i_{p-1}\leq n$. We obtain
$$
\begin{aligned}
&\sum_{1\leq i_1<\cdots <i_{p-1}\leq n} \  
   \sum_{0\,\leq\, a_{i_j}<\,\ell^j_{i_j}}
  \underline{v}^{\underline{Q}(\underline{\ell}^2,\cdots,\underline{\ell}^p;
   i_1,\cdots,i_{p-1})+\underline{\ell}^1}\ot \underline{v}^{\underline{\hat{Q}}}
  \ot v_{i_1}\wedge\cdots\wedge v_{i_{p-1}}\\
&+\sum_{m=1}^{p-1}(-1)^{m}
   \sum_{1\leq i_1<\cdots <i_{p-1}\leq n} \  \ \ 
   \sum_{\substack{ 0\,\leq\, a_{i_j}<\,\ell^{j}_{i_j}\ (
{\mbox{\scriptsize{for} }}
j\in \{1,\ldots, m-1\})\\
\rule[-1.5ex]{0ex}{3.5ex} %strut
0\,\leq\, a_{i_{m}}\leq\, \ell_{i_m}^m+\ell_{i_m}^{m+1}-1\\
     0\,\leq\, a_{i_j}<\,\ell^{j+1}_{i_j} (
{\mbox{\scriptsize{for} }}
j\in \{m+1,\ldots,p\})}}
\\
&\hspace{1in}
\rule[-2ex]{0ex}{6ex} %strut
  \underline{v}^{\underline{Q}(\underline{\ell}^1,\ldots,\underline{\ell}^{m-1},
\underline{\ell}^m+\underline{\ell}^{m+1},\underline{\ell}^{m+2},\ldots,
\underline{\ell}^p; i_1,\ldots,i_{p-1})}\ot
\underline{v}^{\underline{\hat{Q}}}\ot v_{i_1}\wedge
\cdots\wedge v_{i_{p-1}}\\
&+(-1)^{p} \!
\sum_{1\leq i_1<\cdots <i_{p-1}\leq n}\ \ 
   \sum_{0\,\leq\, a_{i_j}<\,\ell^j_{i_j}} \! 
  \ \underline{v}^{\underline{Q}(\underline{\ell}^1,\cdots,\underline{\ell}^{p-1};
   i_2,\cdots,i_{p-1})}\ot \underline{v}^{\underline{\hat{Q}} +\underline{\ell}^p}
  \ot v_{i_1}\wedge\cdots\wedge v_{i_{p-1}}\\
&=\  \Psi_{p-1}(\delta_p(1\ot\underline{v}^{\underline{\ell}^1}\ot
  \cdots\ot \underline{v}^{\underline{\ell}^p}\ot 1)).
\end{aligned}
$$
This finishes the proof of Theorem~\ref{psib-chain-map}.

%%%%%%%%%%%%%%%%%%%%%%%%%%%%%%%%%%%%%%%%%%%%%%%%%%%%%%%%%%%%%%
%%%%%%%%%%%%%%%%%%%%%%%%%%%%%%%%%%%%%%%%%%%%%%%%%%%%%%%%%%%%%%
%%%%%%%%%%%%%%%%%%%%%%%%%%%%%%%%%%%%%%%%%%%%%%%%%%%%%%%%%%%%%%
%%%%%%%%%%%%%%%%%%%%%%%%%%%%%%%%%%%%%%%%%%%%%%%%%%%%%%%%%%%%%%
\quad

%%%%%%%%%%%%%%%%%%%%%%%%%%%%%%%%%%%%%%%%%%%%%%%%%%%%%%%%%%%%%%5
{\sc Acknowledgements.}
The first author thanks Mohamed Barakat for useful conversations.

%%%%%%%%%%%%%%%%%%%%%%%%%%%%%%%%%%%%%%%%%%%%%%%%%%%%%%%%%%%%%%%%%%%%%
%%%%%%%%%%%%%%%%%%%%%%%%%%%%%%%%%%%%%%%%%%%%%%%%%%%%%%%%%%%%%%%%%%%%%

\end{document}